\documentstyle[12pt]{article}

\font\twlgot =eufm10 scaled \magstep1 \font\egtgot =eufm8
\font\sevgot =eufm7 \font\twlmsb =msbm10 scaled \magstep1
\font\egtmsb =msbm8 \font\sevmsb =msbm7

\newfam\gotfam
\def\pgot{\fam\gotfam\twlgot}
\textfont\gotfam\twlgot \scriptfont\gotfam\egtgot
\scriptscriptfont\gotfam\sevgot
\def\got{\protect\pgot}
\newfam\msbfam
\textfont\msbfam\twlmsb \scriptfont\msbfam\egtmsb
\scriptscriptfont\msbfam\sevmsb
\def\Bbb{\protect\pBbb}
\def\pBbb{\relax\ifmmode\expandafter\Bb\else\typeout{You cann't use
Bbb in text mode}\fi}
\def\Bb #1{{\fam\msbfam\relax#1}}

\newcommand{\gd}{{\got d}}
\newcommand{\ccG}{{\got g}}
\newcommand{\gA}{{\got A}}

\def\thebibliography#1{\section*{References}\list
  {[\arabic{enumi}]}{\settowidth\labelwidth{#1}\leftmargin\labelwidth
    \advance\leftmargin\labelsep
    \usecounter{enumi}}
    \def\newblock{\hskip .11em plus .33em minus .07em}
    \sloppy\clubpenalty4000\widowpenalty4000
    \sfcode`\.=1000\relax}

\def\op#1{\mathop{\fam0 #1}\limits}

\newcommand{\beq}{\begin{equation}}
\newcommand{\eeq}{\end{equation}}
\newcommand{\ben}{\begin{eqnarray}}
\newcommand{\een}{\end{eqnarray}}
\newcommand{\be}{\begin{eqnarray*}}
\newcommand{\ee}{\end{eqnarray*}}
\newcommand{\bea}{\begin{eqalph}}
\newcommand{\eea}{\end{eqalph}}
\newcommand{\cA}{{\cal A}}

\newcommand{\cP}{{\cal P}}

\newcommand{\cV}{{\cal V}}
\newcommand{\cC}{{\cal C}}
\newcommand{\cL}{{\cal L}}
\newcommand{\cE}{{\cal E}}

\newcommand{\cS}{{\cal S}}
\newcommand{\cO}{{\cal O}}
\newcommand{\bL}{{\bf L}}

\newcommand{\dl}{\delta}
\newcommand{\la}{\lambda}
\newcommand{\La}{\Lambda}
\newcommand{\f}{\phi}
\newcommand{\om}{\omega}

\newcommand{\m}{\mu}

\newcommand{\G}{\Gamma}
\newcommand{\th}{\theta}

\newcommand{\vt}{\vartheta}
\newcommand{\vf}{\varphi}
\newcommand{\up}{\upsilon}

\newcommand{\di}{{\rm dim\,}}

\newcommand{\si}{\sigma}
\newcommand{\Si}{\Sigma}
\newcommand{\w}{\wedge}

\newcommand{\dr}{\partial}
\newcommand{\ar}{\op\longrightarrow}
\newcommand{\ot}{\otimes}

\newcounter{theorem}
\newcounter{remark}
\newcounter{proposition}
\newcounter{lemma}
\newcounter{corollary}
\newcounter{definition}

\def\theremark{\arabic{remark}}

\def\thedefinition{\arabic{definition}}

\newenvironment{prop}{\refstepcounter{definition} \medskip
\noindent{\it Proposition \thedefinition.}}{\medskip}

\newcommand{\mar}[1]{}

\begin{document}
\hbox{}

{\parindent=0pt

{\large\bf On algebras of gauge transformations in a general
setting}
\bigskip

{\sc G.Sardanashvily}

{\sl Department of Theoretical Physics, Moscow State University,
117234 Moscow, Russia}

\bigskip

{\small {\bf Abstract} We consider a Lagrangian system on a fiber
bundle and its gauge transformations  depending on derivatives of
dynamic variables and gauge parameters of arbitrary order. We say
that gauge transformations form an algebra if they generate a
nilpotent BRST operator. }

 }

\section{Introduction}

In a general setting, one can say that a family of gauge
transformations depending on parameters is an algebra if the Lie
bracket of arbitrary two gauge transformations depending on
different parameter functions is again a gauge transformation
depending on some parameter function. The goal is to formulate
this condition in strict mathematical terms. For instance, gauge
transformations in gauge theory on a principal bundle form a
finite-dimensional real (or complex) Lie algebra. Gauge
transformations of a certain class of field theories constitute
sh-Lie algebras \cite{fulp02}. It may happen that gauge
transformations are not assembled into an algebra or form an
algebra on-shell \cite{gom}.

We consider a Lagrangian system on a smooth fiber bundle $Y\to X$
subject to gauge transformations depending both on derivatives of
dynamic variables of arbitrary order and a finite family of gauge
parameters and their derivatives of arbitrary order.

Let $J^rY$, $r=1,\ldots$, be finite order jet manifolds of
sections of $Y\to X$. In the sequel, the index $r=0$ stands for
$Y$. Given bundle coordinates $(x^\la,y^i)$ on $Y$, jet manifolds
$J^rY$ are endowed with the adapted coordinates
$(x^\la,y^i,y^i_\La)$, where $\La=(\la_k...\la_1)$,
$k=1,\ldots,r$, is a symmetric multi-index. We use the notation
$\la+\La=(\la\la_k...\la_1)$ and
\mar{5.177}\beq
d_\la = \dr_\la + \op\sum_{0\leq|\La|} y^i_{\la+\La}\dr_i^\La,
\qquad d_\La=d_{\la_r}\circ\cdots\circ d_{\la_1}. \label{5.177}
\eeq

In order to describe gauge transformations depending on
parameters, let us consider Lagrangian formalism on the bundle
product
\mar{0681}\beq
E=Y\op \times_X V, \label{0681}
\eeq
where $V\to X$ is a vector bundle whose sections are gauge
parameter functions \cite{noether}. Let $V\to X$ be coordinated by
$(x^\la,\xi^r)$. Then gauge transformations are represented by a
differential operator
\mar{0509}\beq
\up= \op\sum_{0\leq|\La|\leq
m}\up^{i,\La}_r(x^\la,y^i_\Si)\xi^r_\La \dr_i \label{0509}
\eeq
on $E$ (\ref{0681}) which is linear on $V$ and takes its values
into the vertical tangent bundle $VY$ of $Y\to X$. Given a section
$\xi(x)$ of $V\to X$, the pull-back
\mar{0682}\beq
\xi^*\up= \op\sum_{0\leq|\La|\leq m}\up^{i,\La}_r(x^\la,y^i_\Si)
d_\La \xi^r(x) \dr_i \label{0682}
\eeq
of $\up$ (\ref{0509}) onto $Y$ is a gauge transformation depending
on a parameter function $\xi(x)$.

By means of a replacement of even gauge parameters $\xi^r$ and
their jets $\xi^r_\La$ with the odd ghosts $c^r$ and their jets
$c^r_\La$, the operator (\ref{0509}) defines a graded derivation
\mar{0680}\beq
\up= \op\sum_{0\leq|\La|\leq m}\up^{i,\La}_r(x^\la,y^i_\Si)c^r_\La
\dr_i \label{0680}
\eeq
of the algebra of the original even fields and odd ghosts. Its
extension
\mar{0684}\beq
\up= \op\sum_{0\leq|\La|\leq m}\up^{i,\La}_r c^r_\La \dr_i +
u^r\dr_r \label{0684}
\eeq
to ghosts is called the BRST transformation if it is nilpotent.

We say that gauge transformations (\ref{0682}) make up an algebra
if they generate a BRST transformation (\ref{0684}). One can think
of the nilpotency  conditions (\ref{0690}) -- (\ref{0691}) as
being the generalized commutation relations and Jacobi identity,
respectively. This definition is especially convenient for BRST
theory and BV quantization \cite{barn}.

\section{Gauge systems on fiber bundles}

In Lagrangian formalism on a fiber bundle $Y\to X$, Lagrangians
and their Euler--Lagrange operators are represented by elements of
the following graded differential algebra (henceforth GDA).

With the inverse system of jet manifolds
\mar{5.10}\beq
X\op\longleftarrow^\pi Y\op\longleftarrow^{\pi^1_0} J^1Y
\longleftarrow \cdots J^{r-1}Y \op\longleftarrow^{\pi^r_{r-1}}
J^rY\longleftarrow\cdots, \label{5.10}
\eeq
one has the direct system
\mar{5.7}\beq
\cO^*X\op\longrightarrow^{\pi^*} \cO^*Y
\op\longrightarrow^{\pi^1_0{}^*} \cO_1^*Y \ar\cdots \cO^*_{r-1}Y
\op\longrightarrow^{\pi^r_{r-1}{}^*}
 \cO_r^*Y \longrightarrow\cdots \label{5.7}
\eeq
of GDAs $\cO_r^*Y$ of exterior forms on jet manifolds $J^rY$ with
respect to the pull-back monomorphisms $\pi^r_{r-1}{}^*$. Its
direct limit
 $\cO_\infty^*Y$ is a GDA
consisting of all exterior forms on finite order jet manifolds
modulo the pull-back identification.

The projective limit $(J^\infty Y, \pi^\infty_r:J^\infty Y\to
J^rY)$ of the inverse system (\ref{5.10}) is a Fr\'echet manifold.
A bundle atlas $\{(U_Y;x^\la,y^i)\}$ of $Y\to X$ yields the
coordinate atlas
\mar{jet1}\beq
\{((\pi^\infty_0)^{-1}(U_Y); x^\la, y^i_\La)\}, \qquad
{y'}^i_{\la+\La}=\frac{\dr x^\m}{\dr x'^\la}d_\m y'^i_\La, \qquad
0\leq|\La|, \label{jet1}
\eeq
of $J^\infty Y$, where $d_\m$ are the total derivatives
(\ref{5.177}).  Then $\cO^*_\infty Y$ can be written in a
coordinate form where the horizontal one-forms $\{dx^\la\}$ and
the contact one-forms $\{\th^i_\La=dy^i_\La
-y^i_{\la+\La}dx^\la\}$ are generating elements of the
$\cO^0_\infty U_Y$-algebra $\cO^*_\infty U_Y$.

There is the canonical decomposition $\cO^*_\infty
Y=\oplus\cO^{k,m}_\infty Y$ of $\cO^*_\infty Y$ into $\cO^0_\infty
Y$-modules $\cO^{k,m}_\infty Y$ of $k$-contact and $m$-horizontal
forms together with the corresponding projectors $h_k:\cO^*_\infty
Y\to \cO^{k,*}_\infty Y$ and $h^m:\cO^*_\infty Y\to
\cO^{*,m}_\infty Y$. Accordingly, the exterior differential on
$\cO_\infty^* Y$ is split into the sum $d=d_H+d_V$ of the
nilpotent total and vertical differentials
\be
d_H(\f)= dx^\la\w d_\la(\f), \qquad d_V(\f)=\th^i_\La \w
\dr^\La_i\f, \qquad \f\in\cO^*_\infty Y.
\ee

Any finite order Lagrangian
\mar{0512}\beq
L=\cL\om:J^rY\to \op\w^nT^*X, \qquad \om=dx^1\w\cdots\w dx^n,
\qquad n=\di X, \label{0512}
\eeq
is an element of $\cO^{0,n}_\infty Y$, while
\mar{0513}\beq
\dl L=\cE_i\th^i\w\om=\op\sum_{0\leq|\La|}
(-1)^{|\La|}d_\La(\dr^\La_i \cL)\th^i\w\om\in \cO^{1,n}_\infty Y
\label{0513}
\eeq
is its Euler--Lagrange operator taking values into the vector
bundle
\mar{0548}\beq
T^*Y\op\w_Y(\op\w^n T^*X)= V^*Y\op\ot_Y\op\w^n T^*X. \label{0548}
\eeq

A Lagrangian system on a fiber bundle $Y\to X$ is said to be a
gauge theory if its Lagrangian $L$ admits a family of variational
symmetries parameterized by elements of a vector bundle $V\to X$
and its jet manifolds as follows.

Let $\gd\cO^0_\infty Y$ be the $\cO^0_\infty Y$-module of
derivations of the $\Bbb R$-ring $\cO^0_\infty Y$. Any $\vt\in
\gd\cO^0_\infty Y$ yields the graded derivation (the interior
product) $\vt\rfloor\f$ of the GDA $\cO^*_\infty Y$ given by the
relations
\be
&&\vt\rfloor df=\vt(f), \qquad  f\in \cO^0_\infty Y, \\
&& \vt\rfloor(\f\w\si)=(\vt\rfloor \f)\w\si
+(-1)^{|\f|}\f\w(\vt\rfloor\si), \qquad \f,\si\in \cO^*_\infty Y,
\ee
and its derivation (the Lie derivative)
\mar{0515}\ben
&& \bL_\vt\f=\vt\rfloor d\f+ d(\vt\rfloor\f), \qquad \f\in
\cO^*_\infty Y, \label{0515}\\
&& \bL_\vt(\f\w\f')=\bL_\vt(\f)\w\f' +\f\w\bL_\vt(\f'). \nonumber
\een
Relative to an atlas (\ref{jet1}), a derivation
$\vt\in\gd\cO^0_\infty$ reads
\mar{g3}\beq
\vt=\vt^\la \dr_\la + \vt^i\dr_i + \op\sum_{|\La|>0}\vt^i_\La
\dr^\La_i, \label{g3}
\eeq
where $\{\dr_\la,\dr^\La_i\}$ is the dual to the basis $\{dx^\la,
dy^i_\La\}$ with respect to the interior product $\rfloor$
\cite{cmp}.

A derivation $\vt$ is called contact if the Lie derivative
$\bL_\vt$ (\ref{0515}) preserves the contact ideal of the GDA
$\cO^*_\infty Y$ generated by contact forms. A derivation $\up$
(\ref{g3}) is contact iff
\mar{g4}\beq
\vt^i_\La=d_\La(\vt^i-y^i_\m\vt^\m)+y^i_{\m+\La}\vt^\m, \qquad
0<|\La|. \label{g4}
\eeq
Any contact derivation admits the horizontal splitting
\mar{g5}\beq
\vt=\vt_H +\vt_V=\vt^\la d_\la + (\up^i\dr_i + \op\sum_{0<|\La|}
d_\La \up^i\dr_i^\La), \qquad \up^i= \vt^i-y^i_\m\vt^\m.
\label{g5}
\eeq
Its vertical part $\vt_V$  is completely determined by the first
summand
\mar{0641}\beq
\up=\up^i(x^\la,y^i_\La)\dr_i, \qquad 0\leq |\La|\leq k.
\label{0641}
\eeq
This is a section of the pull-back $VY\op\times_Y J^kY\to J^kY$,
i.e., a $k$-order $VY$-valued differential operator on $Y$. One
calls $\up$ (\ref{0641}) a generalized vector field on $Y$.

\begin{prop}  \label{g75} \mar{g75}
The Lie derivative of a Lagrangian $L$ (\ref{0512}) along a
contact derivation $\vt$ (\ref{g5}) fulfills the first variational
formula
\mar{g8'}\beq
\bL_\vt L= \up\rfloor\dl L +d_H(h_0(\vt\rfloor\Xi_L)) +\cL d_V
(\vt_H\rfloor\om), \label{g8'}
\eeq
where $\Xi_L$ is a Lepagean equivalent of $L$ \cite{cmp}
\end{prop}

A contact derivation $\vt$ (\ref{g5}) is called  variational if
the Lie derivative (\ref{g8'}) is $d_H$-exact, i.e., $\bL_\vt
L=d_H\si$, $\si\in \cO^{0,n-1}_\infty$. A glance at the expression
(\ref{g8'}) shows that: (i) $\vt$ (\ref{g5}) is variational only
if it is projected onto $X$; (ii) $\vt$ is variational iff its
vertical part $\vt_V$ is well; (iii) it is variational iff
$\up\rfloor\dl L$ is $d_H$-exact.

By virtue of item (ii), we can restrict our consideration to
vertical contact derivations $\vt=\vt_V$. A generalized vector
field $\up$ (\ref{0641}) is called a variational symmetry of a
Lagrangian $L$ if it generates a variational contact derivation.

Turn now to the notion of a gauge symmetry \cite{noether}. Let us
consider the bundle product $E$ (\ref{0681}) coordinated by
$(x^\la,y^i,\xi^r)$. Given a Lagrangian $L$ on $Y$, let us
consider its pull-back, say again $L$, onto $E$. Let $\vt_E$ be a
contact derivation of the $\Bbb R$-ring $\cO^0_\infty E$, whose
restriction
\mar{0508}\beq
\vt=\vt_E|_{\cO^0_\infty Y}=
\op\sum_{0\leq|\La|}d_\La\up^i\dr_i^\La \label{0508}
\eeq
to $\cO^0_\infty Y\subset \cO^0_\infty E$ is linear in coordinates
$\xi^r_\Xi$. It is determined by a generalized vector field
$\up_E$ on $E$ whose projection
\be
\up:J^kE\ar^{\up_E} VE\to E\op\times_Y VY
\ee
is a linear $VY$-valued differential operator $\up$ (\ref{0509})
on $E$. Let $\vt_E$ be a variational symmetry of a Lagrangian $L$
on $E$, i.e.,
\mar{0552}\beq
\up_E\rfloor \dl L=\up\rfloor \dl L=d_H\si. \label{0552}
\eeq
Then one says that $\up$ (\ref{0509}) is a gauge symmetry of a
Lagrangian $L$.

Note that any differential operator $\up$ (\ref{0509}) defines a
generalized vector field $\up_E=\up$ on $E$ which lives in $VY$
and, consequently, generates a contact derivation $\vt_E=\vt$
(\ref{0508}).

\section{Graded Lagrangian systems}

In order to introduce a BRST operator, let us consider Lagrangian
systems of even and odd variables. We describe odd variables and
their jets on a smooth manifold $X$ as generating elements of the
structure ring of a graded manifold whose body is $X$
\cite{cmp,ijmp}. This definition reproduces the heuristic notion
of jets of ghosts in the field-antifield BRST theory
\cite{barn,bran01}.

Recall that any graded manifold $(\gA,X)$ with a body $X$ is
isomorphic to the one whose structure sheaf $\gA_Q$ is formed by
germs of sections of the exterior product
\mar{g80}\beq
\w Q^*=\Bbb R\op\oplus_X Q^*\op\oplus_X\op\w^2
Q^*\op\oplus_X\cdots, \label{g80}
\eeq
where $Q^*$ is the dual of some real vector bundle $Q\to X$ of
fiber dimension $m$. In field models, a vector bundle $Q$ is
usually given from the beginning. Therefore, we consider graded
manifolds $(X,\gA_Q)$ where the above mentioned isomorphism holds,
and call $(X,\gA_Q)$ the simple graded manifold constructed from
$Q$. The structure ring $\cA_Q$ of sections of $\gA_Q$ consists of
sections of the exterior bundle (\ref{g80}) called graded
functions. Given bundle coordinates $(x^\la,q^a)$ on $Q$ with
transition functions $q'^a=\rho^a_b q^b$, let $\{c^a\}$ be the
corresponding fiber bases for $Q^*\to X$, together with transition
functions $c'^a=\rho^a_bc^b$. Then $(x^\la, c^a)$ is called the
local basis for the graded manifold $(X,\gA_Q)$. With respect to
this basis, graded functions read
\be
f=\op\sum_{k=0}^m \frac1{k!}f_{a_1\ldots a_k}c^{a_1}\cdots
c^{a_k},
\ee
where $f_{a_1\cdots a_k}$ are local smooth real functions on $X$.

Given a graded manifold $(X,\gA_Q)$, let $\gd\cA_Q$ be the
$\cA_Q$-module of $\Bbb Z_2$-graded derivations of the $\Bbb
Z_2$-graded ring of $\cA_Q$, i.e.,
\be
u(ff')=u(f)f'+(-1)^{[u][f]}fu (f'), u\in\gd\cA_Q, \qquad f,f'\in
\cA_Q,
\ee
where $[.]$ denotes the Grassmann parity. Its elements are called
$\Bbb Z_2$-graded (or, simply, graded) vector fields on
$(X,\gA_Q)$. Due to the canonical splitting $VQ= Q\times Q$, the
vertical tangent bundle $VQ\to Q$ of $Q\to X$ can be provided with
the fiber bases $\{\dr_a\}$ which is the dual of $\{c^a\}$. Then a
graded vector field takes the local form $u= u^\la\dr_\la +
u^a\dr_a$, where $u^\la, u^a$ are local graded functions. It acts
on $\cA_Q$ by the rule
\mar{cmp50'}\beq
u(f_{a\ldots b}c^a\cdots c^b)=u^\la\dr_\la(f_{a\ldots b})c^a\cdots
c^b +u^d f_{a\ldots b}\dr_d\rfloor (c^a\cdots c^b). \label{cmp50'}
\eeq
This rule implies the corresponding transformation law
\be
u'^\la =u^\la, \qquad u'^a=\rho^a_ju^j +
u^\la\dr_\la(\rho^a_j)c^j.
\ee
Then one can show \cite{book00,ijmp} that graded vector fields on
a simple graded manifold can be represented by sections of the
vector bundle $\cV_Q\to X$ which is locally isomorphic to the
vector bundle $\w Q^*\ot_X(Q\oplus_X TX)$.

Using this fact, we can introduce graded exterior forms on the
graded manifold $(X,\gA_Q)$ as sections of the exterior bundle
$\op\w\cV^*_Q$, where $\cV^*_Q\to  X$ is the $\w Q^*$-dual of
$\cV_Q$. Relative to the dual local bases $\{dx^\la\}$ for $T^*X$
and $\{dc^b\}$ for $Q^*$, graded one-forms read
\be
\f=\f_\la dx^\la + \f_adc^a,\qquad \f'_a=\rho^{-1}{}_a^b\f_b,
\qquad \f'_\la=\f_\la +\rho^{-1}{}_a^b\dr_\la(\rho^a_j)\f_bc^j.
\ee
The duality morphism is given by the interior product
\be
u\rfloor \f=u^\la\f_\la + (-1)^{[\f_a]}u^a\f_a.
\ee
Graded exterior forms constitute the bigraded differential algebra
(henceforth BGDA) $\cC^*_Q$ with respect to the bigraded exterior
product $\w$ and the exterior differential $d$. The standard
formulae of a BGDA hold.

Since the jet bundle $J^rQ\to X$ of a vector bundle $Q\to X$ is a
vector bundle, let us consider the simple graded manifold
$(X,\gA_{J^rQ})$ constructed from $J^rQ\to X$. Its local basis is
$\{x^\la,c^a_\La\}$, $0\leq |\La|\leq r$, together with the
transition functions
\mar{+471}\beq
c'^a_{\la +\La}=d_\la(\rho^a_j c^j_\La), \qquad d_\la=\dr_\la +
\op\sum_{|\La|<r}c^a_{\la+\La} \dr_a^\La, \label{+471}
\eeq
where $\dr_a^\La$ are the duals of $c^a_\La$. Let $\cC^*_{J^rQ}$
be the BGDA of graded exterior forms on the graded manifold
$(X,\gA_{J^rQ})$. A linear bundle morphism $\pi^r_{r-1}:J^rQ \to
J^{r-1}Q$ yields the corresponding monomorphism of BGDAs
$\cC^*_{J^{r-1}Q}\to \cC^*_{J^rQ}$. Hence, there is the direct
system of BGDAs
\mar{g205}\beq
\cC^*_Q\ar^{\pi^{1*}_0} \cC^*_{J^1Q}\cdots
\ar^{\pi^r_{r-1}{}^*}\cC^*_{J^rQ}\ar\cdots. \label{g205}
\eeq
Its direct limit $\cC^*_\infty Q$ consists of graded exterior
forms on graded manifolds $(X,\gA_{J^rQ})$, $r\in\Bbb N$, modulo
the pull-back identification, and it inherits the BGDA operations
intertwined by the monomorphisms $\pi^r_{r-1}{}^*$. It is a
$C^\infty(X)$-algebra locally generated by the elements $(1,
c^a_\La, dx^\la,\th^a_\La=dc^a_\La -c^a_{\la +\La}dx^\la)$,
$0\leq|\La|$.

 In order to regard even and odd dynamic variables on the
same footing,  let $Y\to X$ be hereafter an affine bundle, and let
$\cP^*_\infty Y\subset \cO^*_\infty Y$ be the
$C^\infty(X)$-subalgebra of exterior forms whose coefficients are
polynomial in the fiber coordinates $y^i_\La$ on jet bundles $J^r
Y\to X$. Let us consider the product
\mar{0670}\beq
\cS^*_\infty=\cC_\infty^*Q\w\cP^*_\infty Y \label{0670}
\eeq
of graded algebras $\cC_\infty^*Q$ and $\cP^*_\infty Y$ over their
common graded subalgebra $\cO^*X$ of exterior forms on $X$
\cite{cmp}. It consists of the elements
\be
\op\sum_i \psi_i\ot\f_i, \qquad \op\sum_i \f_i\ot\psi_i, \qquad
\psi\in \cC^*_\infty Q, \qquad \f\in \cP^*_\infty Y,
\ee
modulo the commutation relations
\mar{0442}\ben
&&\psi\ot\f=(-1)^{|\psi||\f|}\f\ot\psi, \qquad
\psi\in \cC^*_\infty Q, \qquad \f\in \cP^*_\infty Y, \label{0442}\\
&& (\psi\w\si)\ot\f=\psi\ot(\si\w\f), \qquad \si\in \cO^*X.
\nonumber
\een
They are  endowed with the total form degree $|\psi|+|\f|$ and the
total Grassmann parity $[\psi]$. Their multiplication
\mar{0440}\beq
(\psi\ot\f)\w(\psi'\ot\f'):=(-1)^{|\psi'||\f|}(\psi\w\psi')\ot
(\f\w\f'). \label{0440}
\eeq
obeys the relation
\be
\vf\w\vf' =(-1)^{|\vf||\vf'| +[\vf][\vf']}\vf'\w \vf, \qquad
\vf,\vf'\in \cS^*_\infty,
\ee
and makes $\cS^*_\infty$ (\ref{0670}) into a bigraded $C^\infty
(X)$-algebra.  For instance, elements of the ring $S^0_\infty$ are
polynomials of $c^a_\La$ and $y^i_\La$ with coefficients in
$C^\infty(X)$.

The algebra $\cS^*_\infty$ is provided with the exterior
differential
\mar{0441}\beq
d(\psi\ot\f):=(d_\cC\psi)\ot\f +(-1)^{|\psi|}\psi\ot(d_\cP\f),
\qquad \psi\in \cC^*_\infty, \qquad \f\in \cP^*_\infty,
\label{0441}
\eeq
where $d_\cC$ and $d_\cP$ are exterior differentials on the
differential algebras $\cC^*_\infty Q$ and $\cP^*_\infty Y$,
respectively. It obeys the relations
\be
 d(\vf\w\vf')= d\vf\w\vf' +(-1)^{|\vf|}\vf\w d\vf', \qquad
\vf,\vf'\in \cS^*_\infty,
\ee
and makes $\cS^*_\infty$ into a BGDA, which is locally generated
by the elements
\be
(1, c^a_\La, y^i_\La,
dx^\la,\th^a_\La=dc^a_\La-c^a_{\la+\La}dx^\la,\th^i_\La=
dy^i_\La-y^i_{\la+\La}dx^\la), \qquad 0\leq |\La|.
\ee

Hereafter, let the collective symbols $s^A_\La$ and $\th^A_\La$
stand both for even and odd generating elements $c^a_\La$,
$y^i_\La$, $\th^a_\La$, $\th^i_\La$ of the $C^\infty(X)$-algebra
$\cS^*_\infty$ which, thus, is locally generated by $(1,s^A_\La,
dx^\la, \th^A_\La)$, $|\La|\geq 0$. We agree to call elements of
$\cS^*_\infty$ the graded exterior forms on $X$.

Similarly to $\cO^*_\infty Y$, the BGDA $\cS^*_\infty$ is
decomposed into $\cS^0_\infty$-modules $\cS^{k,r}_\infty$ of
$k$-contact and $r$-horizontal graded forms together with the
corresponding projections $h_k$ and $h^r$. Accordingly, the
exterior differential $d$ (\ref{0441}) on $\cS^*_\infty$ is split
into the sum $d=d_H+d_V$ of the total and vertical differentials
\be
d_H(\f)=dx^\la\w d_\la(\f), \qquad d_V(\f)=\th^A_\La\w\dr^\La_A
\f, \qquad \f\in \cS^*_\infty.
\ee
One can think of the elements
\be
L=\cL\om\in \cS^{0,n}_\infty, \qquad \dl (L)= \op\sum_{|\La|\geq
0}
 (-1)^{|\La|}\th^A\w d_\La (\dr^\La_A L)\in \cS^{0,n}_\infty
\ee
as being a graded Lagrangian and its Euler--Lagrange operator,
respectively.

\section{BRST symmetry}

A graded derivation  $\vt\in\gd \cS^0_\infty$ of the $\Bbb R$-ring
$\cS^0_\infty$ is said to be contact if the Lie derivative
$\bL_\vt$ preserves the ideal of contact graded forms of the BGDA
$\cS^*_\infty$. With respect to the local basis $(x^\la,s^A_\La,
dx^\la,\th^A_\La)$ for the BGDA $\cS^*_\infty$, any contact graded
derivation takes the form
\mar{g105}\beq
\vt=\vt_H+\vt_V=\vt^\la d_\la + (\vt^A\dr_A +\op\sum_{|\La|>0}
d_\La\vt^A\dr_A^\La), \label{g105}
\eeq
where $\vt^\la$, $\vt^A$ are local graded functions \cite{cmp}.
The interior product $\vt\rfloor\f$ and the Lie derivative
$\bL_\vt\f$, $\f\in\cS^*_\infty$, are defined by the same formulae
\be
&& \vt\rfloor \f=\vt^\la\f_\la + (-1)^{[\f_A]}\vt^A\f_A, \qquad
\f\in \cS^1_\infty,\\
&& \vt\rfloor(\f\w\si)=(\vt\rfloor \f)\w\si
+(-1)^{|\f|+[\f][\vt]}\f\w(\vt\rfloor\si), \qquad \f,\si\in
\cS^*_\infty \\
&& \bL_\vt\f=\vt\rfloor d\f+ d(\vt\rfloor\f), \qquad
\bL_\vt(\f\w\si)=\bL_\vt(\f)\w\si
+(-1)^{[\vt][\f]}\f\w\bL_\vt(\si).
\ee
as those on a graded manifold. One can justify that any vertical
contact graded derivation $\vt$ (\ref{g105}) satisfies the
relations
\mar{g232}\beq
\vt\rfloor d_H\f=-d_H(\vt\rfloor\f), \qquad
\bL_\vt(d_H\f)=d_H(\bL_\vt\f), \qquad \f\in\cS^*_\infty.
\label{g232}
\eeq

\begin{prop}  \label{g106} \mar{g106}
The Lie derivative $\bL_\vt L$ of a Lagrangian $L$ along a contact
graded derivation $\vt$ (\ref{g105}) fulfills the first
variational formula
\mar{g107}\beq
\bL_\vt L= \vt_V\rfloor\dl L +d_H(h_0(\vt\rfloor \Xi_L)) + d_V
(\vt_H\rfloor\om)\cL, \label{g107}
\eeq
where $\Xi_L=\Xi+L$ is a Lepagean equivalent of a graded
Lagrangian $L$ \cite{cmp}.
\end{prop}

A contact graded derivation $\vt$ is said to be variational if the
Lie derivative (\ref{g107}) is $d_H$-exact. A glance at the
expression (\ref{g107}) shows that: (i) A contact graded
derivation $\vt$ is variational only if it is projected onto $X$,
and (ii) $\vt$ is variational iff its vertical part $\vt_V$ is
well. Therefore, we restrict our consideration to vertical contact
graded derivations
\mar{0672}\beq
\vt=\op\sum_{0\leq|\La|} d_\La\up^A\dr_A^\La. \label{0672}
\eeq
Such a derivation is completely defined by its first summand
\mar{0673}\beq
\up=\up^A(x^\la,s^A_\La)\dr_A, \qquad 0\leq|\La|\leq k,
\label{0673}
\eeq
which is also a graded derivation of $\cS^0_\infty$. It is called
the generalized graded vector field. A glance at the first
variational formula (\ref{g107}) shows that $\vt$ (\ref{0672}) is
variational iff $\up\rfloor \dl L$ is $d_H$-exact.

A vertical contact graded derivation $\vt$ (\ref{0672}) is said to
be nilpotent if
\mar{g133}\beq
\bL_\up(\bL_\up\f)= \op\sum_{|\Si|\geq 0,|\La|\geq 0 }
(\up^B_\Si\dr^\Si_B(\up^A_\La)\dr^\La_A +
(-1)^{[s^B][\up^A]}\up^B_\Si\up^A_\La\dr^\Si_B \dr^\La_A)\f=0
\label{g133}
\eeq
for any horizontal graded form $\f\in S^{0,*}_\infty$ or,
equivalently, $\vt\circ\vt)(f)=0$ for any graded function $f\in
\cS^0_\infty$. One can show that $\vt$ is nilpotent only if it is
odd and iff the equality
\mar{0688}\beq
\vt(\up^A)=\op\sum_{|\Si|\geq 0}
\up^B_\Si\dr^\Si_B(\up^A)=0 \label{0688}
\eeq
holds for all $\up^A$ \cite{cmp}.

Return now to the original gauge system on a fiber bundle $Y$ with
a Lagrangian $L$ (\ref{0512}) and a gauge symmetry $\up$
(\ref{0509}). For the sake of simplicity, $Y\to X$ is assumed to
be affine. Let us consider the BGDA $\cS^*_\infty=\cC^*_\infty
V\w\cP^*_\infty Y$ locally generated by $(1,c^r_\La, dx^\la,
y^i_\La, \th^r_\La, \th^i_\La)$. Let $L\in \cO^{0,n}_\infty Y$ be
a polynomial in $y^i_\La$, $0\leq |L|$. Then it is a graded
Lagrangian $L\in \cP^{0,n}_\infty Y\subset \cS^{0,n}_\infty$ in
$\cS^*_\infty$. Its gauge symmetry $\up$ (\ref{0509}) gives rise
to the generalized vector field $\up_E=\up$ on $E$, and the latter
defines the generalized graded vector field $\up$ (\ref{0673}) by
the formula (\ref{0680}). It is easily justified that the contact
graded derivation $\vt$ (\ref{0672}) generated by $\up$
(\ref{0680}) is variational for $L$. It is odd, but need not be
nilpotent. However, one can try to find a nilpotent contact graded
derivation (\ref{0672}) generated by some generalized graded
vector field (\ref{0684}) which coincides with $\vt$ on
$\cP^*_\infty Y$. We agree to call it the BRST operator.

In this case, the nilpotency conditions (\ref{0688}) read
\mar{0690,1}\ben
&& \op\sum_\Si d_\Si(\op\sum_\Xi\up^{i,\Xi}_rc^r_\Xi)
\op\sum_\La\dr^\Si_i (\up^{j,\La}_s)c^s_\La +\op\sum_\La d_\La
(u^r)\up^{j,\La}_r
=0, \label{0690}\\
&& \op\sum_\La(\op\sum_\Xi d_\La(\up^{i,\Xi}_r c^r_\Xi)\dr^\La_i
+d_\La (u^r)\dr_r^\La)u^q=0\label{0691}
\een
for all indices $j$ and $q$. They are equations for graded
functions $u^r\in\cS^0_\infty$. Since these functions are
polynomials
\mar{0693}\beq
u^r=u_{(0)}^r + \op\sum_\G u_{(1)p}^{r,\G} c^p_\G +
\op\sum_{\G_1,\G_2} u_{(2)p_1p_2}^{r,\G_1\G_2}
c^{p_1}_{\G_1}c^{p_2}_{\G_2} +\cdots \label{0693}
\eeq
in $c^s_\La$, the equations (\ref{0690}) -- (\ref{0691}) take the
form
\mar{0694}\ben
&& \op\sum_\Si d_\Si(\op\sum_\Xi\up^{i,\Xi}_rc^r_\Xi)
\op\sum_\La\dr^\Si_i (\up^{j,\La}_s)c^s_\La +\op\sum_\La d_\La
(u_{(2)}^r)\up^{j,\La}_r
=0, \label{0694a}\\
&& \op\sum_\La d_\La (u_{(k\neq 2)}^r)\up^{j,\La}_r =0,
\label{0694b}\\
&& \op\sum_\La\op\sum_\Xi d_\La(\up^{i,\Xi}_r c^r_\Xi)\dr^\La_i
u_{(k-1)}^q +\op\sum_{m+n-1=k}d_\La (u_{(m)}^r)\dr_r^\La u_{(n)}^q
=0. \label{0694c}
\een

If the equations (\ref{0694a}) -- (\ref{0694c}) have a solution,
i.e., the (nilpotent) BRST operator exists, one can think of the
equalities (\ref{0694a}) and (\ref{0694c}) (and, consequently, the
nilpotency conditions (\ref{0690}) -- (\ref{0691})) as being the
generalized commutation relations and generalized Jacobi
identities of gauge transformations, respectively.

Indeed, the relation (\ref{0694a}) for components $\up^i_r$ takes
the form of the familiar Lie bracket
\be
 \op\sum_\Si
[d_\Si(\up^i_p)\dr^\Si_i\up^j_q-
d_\Si(\up^i_q)\dr^\Si_i\up^j_p]=-2u^r_{(2)pq}\up^j_r,
\ee
where $-2u^r_{(2)pq}$ are generalized structure constants
depending on dynamic variables $y^i$ and their jets $y^i_\La$ in
general. For instance, let us assume that all $\up^i_r$ are linear
in $y^i_\La$. Then $u^r_{(2)pq}$ are independent of these
variables. Let $u^r_{(m\neq 2)}=0$. In this case, the relation
(\ref{0694c}) reduces to the familiar Jacobi identity
\be
u^r_{(2)pq}u^j_{(2)rs} + u^r_{(2)qs}u^j_{(2)rp}
+u^r_{(2)sp}u^j_{(2)rq}=0.
\ee

Let us note that any Lagrangian $L$ have gauge symmetries. In
particular, there always exist trivial gauge symmetries
\be
\up=\op\sum_\La T_r^{j,i,\La}\cE_j\xi^r_\La\dr_i, \qquad
T_r^{j,i,\La}=-T_r^{i,j,\La},
\ee
vanishing on-shell. In a general setting, one therefore can
require that the nilpotency conditions (\ref{0690}) --
(\ref{0691}) hold on-shell, i.e., gauge transformations form an
algebra on-shell.

\section{Example}

Let us consider the gauge theory of principal connections on a
principal bundle $P\to X$ with a structure Lie group $G$. These
connections are represented by sections of the quotient
\mar{0654}\beq
C=J^1P/G\to X. \label{0654}
\eeq
This is an affine bundle coordinated by $(x^\la, a^r_\la)$ such
that, given a section $A$ of $C\to X$, its components
$A^r_\la=a^r_\la\circ A$ are coefficients of the familiar local
connection form (i.e., gauge potentials). Let $J^\infty C$ be the
infinite order jet manifold of $C\to X$ coordinated by
$(x^\la,a^r_{\La\la})$, $0\leq |\La|$. We consider the GDA
$\cO^*_\infty C$.

Infinitesimal generators of one-parameter groups of automorphisms
of a principal bundle $P$ are $G$-invariant projectable vector
fields on $P\to X$. They are associated to sections of the vector
bundle $T_GP=TP/G\to X$. This bundle is provided with the
coordinates $(x^\la,\dot x^\la,\xi^r)$ with respect to the fibre
bases $\{\dr_\la, e_r\}$ for $T_GP$, where $\{e_r\}$ is the basis
for the right Lie algebra $\ccG$ of $G$ such that
$[e_p,e_q]=c^r_{pq}e_r.$ If
\mar{0652}\beq
 u=u^\la\dr_\la
+u^r e_r, \qquad v=v^\la\dr_\la +v^r e_r, \label{0652}
\eeq
are sections of $T_GP\to X$, their bracket reads
\mar{0651}\beq
[u,v]=(u^\m\dr_\m v^\la -v^\m\dr_\m u^\la)\dr_\la +(u^\la\dr_\la
v^r - v^\la\dr_\la u^r +c^r_{pq}u^pv^q)e_r. \label{0651}
\eeq
Any section $u$ of the vector bundle $T_GP\to X$ yields the vector
field
\mar{0653}\beq
u_C=u^\la\dr_\la +(c^r_{pq}a^p_\la u^q +\dr_\la u^r- a^r_\m\dr_\la
u^\m)\dr^\la_r \label{0653}
\eeq
on the bundle of principal connections $C$ (\ref{0654}). It is an
infinitesimal generator of a one-parameter group of automorphisms
of $C$ \cite{book00}. Let us consider the bundle product
\mar{0659}\beq
E=C\op\times_X T_GP, \label{0659}
\eeq
coordinated by $(x^\la, \tau^\la=\dot x^\la, \xi^r, a^r_\la)$. It
can be provided with the generalized vector field
\mar{0660}\beq
\up_E= \up=(c^r_{pq}a^p_\la \xi^q + \xi^r_\la
-a^r_\m\tau^\m_\la-\tau^\m a_{\m\la}^r)\dr^\la_r. \label{0660}
\eeq

Following the procedure in Sections 3 -- 4, we replace parameters
$\xi^r$ and $\tau^\la$ with the odd ghosts $c^r$ and $c^\la$,
respectively, and obtain the generalized graded vector field
\mar{0695}\beq
\up= (c^r_{pq}a^p_\la c^q + c^r_\la -a^r_\m c^\m_\la-c^\m
a_{\m\la}^r)\dr^\la_r +(-\frac12c^r_{pq}c^pc^q -c^\m c^r_\m)\dr_r
+ c^\la_\m c^\m\dr_\la \label{0695}
\eeq
such that the vertical contact graded derivations (\ref{0672})
generated by $\up$ (\ref{0695}) is nilpotent, i.e., it is a BRST
operator.


\begin{thebibliography}{ddd}

\bibitem{barn} G.Barnich, F.Brandt  and M.Henneaux, Local
BRST cohomology in gauge theories, {\it Phys. Rep.} {\bf 338}, 439
(2000).

\bibitem{noether} D.Bashkirov, G.Giachetta, L.Mangiarotti and
G.Sardanashvily, Noether's second theorem in a general setting.
Reducible gauge theories, {\it E-print arXiv}: math.DG/0411070.

\bibitem{bran01} F.Brandt, Jet coordinates for local BRST
cohomology, {\it Lett. Math. Phys.} {\bf 55}, 149 (2001).

\bibitem{fulp02} R.Fulp, T.Lada and J.Stasheff, Sh-Lie
algebras induced by gauge transformations, {\it Comm. Math. Phys.}
{\bf 231}, 25 (2002); {\it E-print arXiv}: math.QA/0012106.

\bibitem{cmp} G.Giachetta, L.Mangiarotti and
G.Sardanashvily, Lagrangian supersymmetries depending on
derivatives. Global analysis and cohomology, {\it Commun. Math.
Phys.} (accepted); {\it E-print arXiv}: hep-th/0407185.

\bibitem{gom} J.Gomis, J.Par\'\i s, J. and S.Samuel,
Antibracket, antifields and gauge theory quantization, {\it Phys.
Rep} {\bf 295}, 1 (1995).

\bibitem{book00} L.Mangiarotti and G.Sardanashvily, {\it
Connections in Classical and Quantum Field Theory} (World
Scientific, Singapore, 2000).

\bibitem{ijmp} G.Sardanashvily, SUSY-extended field theory,
{\it Int. J. Mod. Phys. A} {\bf 15}, 3095 (2000).

\end{thebibliography}
\end{document}